\documentclass[english]{IEEEtran}
\usepackage[T1]{fontenc}
\usepackage[latin9]{inputenc}
\usepackage{babel}

\usepackage[top=0.75in, bottom=1in, left=0.625in, right=0.625in]{geometry}

\usepackage{epsfig,psfrag}
\usepackage{pst-all}
\usepackage{amsmath,amsthm,amssymb,amsfonts,upref,cite,epsf,color,bm}
\usepackage{graphicx}
\usepackage{color}

\newcommand{\eq}{\,=\,}

\newcommand{\be}{\begin{equation}}
\newcommand{\ee}{\end{equation}}
\newcommand{\ist}{\hspace*{.2mm}}
\newcommand{\rmv}{\hspace*{-.2mm}}
\newcommand{\estmsex}{\varepsilon(\hat{\mathbf{x}}(\cdot);\mathbf{x}) }
\newcommand{\estbiasx}{\mathbf{b}(\hat{\mathbf{x}}(\cdot);\mathbf{x}) }
\newcommand{\estvarx}{v(\hat{\mathbf{x}}(\cdot);\mathbf{x}) }

\newcommand{\RKHSc}{ \mathcal{H}(R_{\mathbf{x}_{0}}^{\mathcal{K}}) }

\newcommand{\x}{{\bf x}}
\newcommand{\xtrue}{\x}
\newcommand{\y}{{\bf y}}
\newcommand{\normconstgauss}{\frac{1}{(2 \pi \sigma^{2})^{M/2}} }
\DeclareMathOperator*{\argmax}{arg\,max}
\DeclareMathOperator{\supp}{supp}

\newcommand{\CRBfull} {{C}ram\'{e}r--{R}ao bound}

\definecolor{hl}{rgb}{0.7,0,0}
\definecolor{fr}{rgb}{0.9,0.5,0}
\definecolor{ret}{rgb}{0,.5,0}

\def\ML_est{\hat{\mathbf{x}}_{\text{ML}}}
\allowdisplaybreaks

\begin{document}

\title{A Lower Bound on the Estimator Variance\\[-.5mm]
for the Sparse Linear Model\vspace{2mm}}

\author{\emph{Sebastian Schmutzhard}$^{\ist 1}\rmv$, 
\emph{Alexander Jung}$^2\rmv$, 
\emph{Franz Hlawatsch}$^2\rmv$, 
\emph{Zvika Ben-Haim}$^{\ist 3}\rmv$, and 
\emph{Yonina C.\ Eldar}$^3$\\[3.5mm]
\normalsize $^1$\normalsize NuHAG, Faculty of Mathematics, University of Vienna\\[-.4mm]
\normalsize 
%% Alserbachstrasse 23/4, 
A-1090 Vienna, Austria; e-mail: sebastian.schmutzhard@univie.ac.at\\[.8mm]
$^2$Institute of Communications and Radio-Frequency Engineering, Vienna University of Technology\\[-.4mm]
\normalsize 
%% Gusshausstrasse 25/389, 
A-1040 Vienna, Austria; e-mail: \{ajung, fhlawats\}@nt.tuwien.ac.at\\[.8mm]
$^3$\normalsize Technion---Israel Institute of Technology\\[-.5mm]
Haifa 32000, Israel; e-mail: \{zvikabh@tx, yonina@ee\}.technion.ac.il
\thanks{This work was supported by the FWF under Grants S10602-N13 (Signal and Information Representation) and S10603-N13 (Statistical Inference) within the National Research Network SISE, by the Israel Science Foundation under Grant 1081/07, and by the European Commission under the FP7 Network of Excellence in Wireless COMmunications NEWCOM++ 
(contract no. 216715).}
}

\maketitle

\thispagestyle{empty}
\pagestyle{empty}

\renewcommand{\baselinestretch}{0.935}\small\normalsize

\begin{abstract}
We study the performance of estimators of a sparse nonrandom vector based on an observation which is linearly transformed and corrupted by 
additive 
white Gaussian noise. 
Using the \emph{reproducing kernel Hilbert space} framework, we derive 
a new lower bound on the estimator variance for a given differentiable bias function (including the unbiased case) and an almost arbitrary transformation matrix (including the
underdetermined case considered in compressed sensing theory). For the special case of a sparse vector corrupted by white Gaussian noise---i.e., without a linear transformation---and 
unbiased estimation, our lower bound improves on previously proposed bounds. 
\end{abstract}

\begin{keywords}
Sparsity, parameter estimation, sparse linear model, denoising, variance bound, reproducing kernel Hilbert space, RKHS.
\vspace{-3mm}
\end{keywords}

%% \pagebreak %%%%%%%%%%

%%%%%%%%%%%%%%%%%%%%%%%%%%%%%%%%%%%
\section{Introduction}
%%%%%%%%%%%%%%%%%%%%%%%%%%%%%%%%%%%

\vspace{.5mm}

We study the problem of estimating a nonrandom parameter vector $\mathbf{x} \!\in\! \mathbb{R}^{N}\!$
which is 
%% known to be 
sparse, i.e., at most $S$ of its entries are nonzero, where $1 \!\leq\! S \!<\! N$ (typically $S \!\ll\! N$). We thus 
\vspace{.5mm}
have
\be
\label{equ_SLM_parameter}
\mathbf{x} \!\in\!  \mathcal{X}_{S} \,, \quad\; \text{with} \;\; \mathcal{X}_{S} \triangleq \big\{\mathbf{x}' \rmv\!\in\! \mathbb{R}^{N} \big| \ist {\|\mathbf{x}' \|}_{0} \rmv\leq\rmv S \big\} \,,
\vspace{.5mm}
\ee 
where ${\| \mathbf{x} \|}_{0}$ denotes the number of nonzero entries of $\mathbf{x}$. 
While the sparsity degree $S$ is assumed to be known, the set of positions of the nonzero entries of $\mathbf{x}$ (denoted by $\rm{supp}(\mathbf{x})$) is unknown. The estimation of $\mathbf{x}$ is based on the observed vector $\mathbf{y} \!\in\! \mathbb{R}^{M}\!$ 
%% which is a linearly transformed and noisy version of $\mathbf{x}$: 
given 
\vspace{-.5mm}
by
\be
\label{equ_observation_model}
\mathbf{y} \ist=\ist \mathbf{H} \mathbf{x} + \mathbf{n} \,,
\vspace{.2mm}
\ee
with a known 
%% deterministic 
system matrix $\mathbf{H} \!\in\! \mathbb{R}^{M \times N}\!$ and white Gaussian noise $\mathbf{n} \sim \mathcal{N}(\mathbf{0}, \sigma^{2} \mathbf{I})$ 
with known variance $\sigma^{2} \!>\! 0$. The matrix $\mathbf{H}$ is arbitrary except that it is assumed to satisfy the standard 
\vspace{-.5mm}
requirement %% BLIBLABLU
\begin{equation} 
\rm{spark}(\mathbf{H}) \rmv>\rmv  S \,, 
\label{equ_spark_cond}
\end{equation} 
where $\rm{spark}(\mathbf{H})$ denotes the minimum number of linearly dependent columns of $\mathbf{H}$ \cite{GreedisGood}.
%(This assumption will be needed in Section \ref{sec:isometry} to relate the given sparse estimation problem to a nonsparse problem.)
The observation 
\pagebreak %%%%%%%%%%
model \eqref{equ_observation_model} together with \eqref{equ_SLM_parameter} will be referred to as the 
\emph{sparse linear model} (SLM).
Note that we also allow $M \!<\! N$ 
%% may be larger than, equal to, or smaller than $N$ 
(this case is relevant to compressed sensing methods \cite{Don06,GreedisGood});
however, condition \eqref{equ_spark_cond} implies that $M \!\ge\! S$.
The case of correlated Gaussian noise $\mathbf{n}$ with a known nonsingular correlation matrix can be reduced to the SLM by means of a noise whitening transformation.
An important special case of the SLM is given by $\mathbf{H} \!=\! \mathbf{I}$ (so that $M \!=\! N$), 
\vspace{-1mm}
i.e., 
\be
\label{equ_observation_model_SSNM}
\mathbf{y} \ist=\ist \mathbf{x} + \mathbf{n} \,,
\ee
where again $\mathbf{x} \!\in\!  \mathcal{X}_{S}$ and $\mathbf{n} \sim \mathcal{N}(\mathbf{0}, \sigma^{2} \mathbf{I})$. This will be referred to as the \emph{sparse signal in noise model} (SSNM).

Lower bounds on the estimation variance for the SLM have been studied previously. 
In particular, the \CRBfull\ (CRB) for the SLM was derived in \cite{ZvikaCRB}.
For the SSNM \eqref{equ_observation_model_SSNM}, lower and upper bounds on the minimum variance of unbiased estimators were derived in \cite{AlexZvikaICASSP}.
A problem with the lower bounds of \cite{ZvikaCRB} and \cite{AlexZvikaICASSP} is the fact that they 
exhibit a discontinuity when passing from the case ${\| \mathbf{x} \|}_{0} \!=\! S$ to the case ${\| \mathbf{x} \|}_{0} \!<\! S$. 

In this paper, we use the mathematical framework of \emph{reproducing kernel Hilbert spaces} (RKHS) \cite{aronszajn1950,Parzen59,Duttweiler73b} 
to derive a novel lower variance bound for the SLM. 
%It is worth noting that some established lower 
%% variance 
%bounds like the CRB and the Hammersley-Chapman-Robbins bound \cite{GormanHero} have pleasing geometric interpretations 
%in the RKHS framework. 
The RKHS framework allows pleasing geometric interpretations of existing bounds, including the CRB, the Hammersley-Chapman-Robbins bound \cite{GormanHero}, 
and the Barankin bound \cite{barankin49}.
The bound we derive here holds for estimators with a given differentiable bias function.
%%  which can be specified arbitrarily. 
For the SSNM, in particular, we obtain a lower bound for unbiased estimators 
%% (i.e., estimators whose bias is identically zero) 
which is tighter
%%  (higher) 
than the bounds in \cite{AlexZvikaICASSP} and, moreover, everywhere continuous. As we will show, RKHS theory relates the 
%% variance 
bound for the SLM to that obtained for the 
%% simpler 
linear model without a sparsity assumption. We note that the RKHS framework has been previously applied to estimation \cite{Parzen59,Duttweiler73b} but, 
to the best of our knowledge, not to the SLM.
%% it has not so far been used for the estimation of sparse vectors. 
%% To the best of our knowledge, RKHS theory has not been previously applied to the estimation of sparse vectors !!!!!.

This paper is organized as follows. 
In Section \ref{sec_toolkit}, we review some fundamentals of parameter estimation.
Relevant elements of RKHS theory are summarized in Section \ref{sec_rkhs}.
In Section \ref{sec_lower_bounds}, we use RKHS theory to derive a lower variance bound for the SLM. 
Section \ref{sec_lower_bounds_SSNM} considers the special case of unbiased estimation within the SSNM. 
Section \ref{sec_numerical} presents a numerical comparison of the new bound with the variance of two established estimation schemes.

\vspace{-.5mm}

%%%%%%%%%%%%%%%%%%%%%%%%%%%%%%%%%%%
\section{Basic Concepts}\label{sec_toolkit}
%%%%%%%%%%%%%%%%%%%%%%%%%%%%%%%%%%%

\vspace{.8mm}

We first review some basic concepts of parameter estimation \cite{kay}. %%% BLIBLABLU
Let $\mathbf{x} \!\in\! \mathcal{X} \!\subseteq\! \mathbb{R}^N\!$ be the nonrandom parameter vector to be estimated, 
%% (we will denote it by $\mathbf{x}_0$ when we wish to emphasize that it is considered fixed) 
$\mathbf{y} \!\in\! \mathbb{R}^M\!$ the observed vector, and $f ( \mathbf{y};  \x)$ the probability density function (pdf) of $\y$, 
parameterized by $\mathbf{x}$. For the SLM, $\mathcal{X} \!=\! \mathcal{X}_S$ as defined in \eqref{equ_SLM_parameter} and
\begin{equation}
\label{equ_pdf}
f ( \mathbf{y};  \x) \,=\, \normconstgauss \, \exp\rmv\rmv \bigg( \!\!-\frac{1}{2 \sigma^{2}} \| \mathbf{y} \!-\! \mathbf{H}\x \|_2^{2} \bigg) \,.
\vspace{-1mm}
\end{equation}

\subsection{Minimum-Variance Estimators}\label{sec:min_var}
%%%%%%%%%%%%%%%%%%%%%%%%%%%%%%%%%%%
%%  with Prescribed Bias

\vspace{.5mm}

%% Consider an estimator $\hat{\mathbf{x}}(\mathbf{y})$ of the parameter vector $\mathbf{x}$.
%%  is a function $\hat{\mathbf{x}}(\cdot) : \mathbb{R}^{M} \!\rmv\rightarrow \mathbb{R}^{N}$
%% that maps the observation $\mathbf{y}$ to an estimate $\hat{\mathbf{x}}$ of $\x$. 
%% Note that even though $\x$ is known to be $S$-sparse, $\hat{\mathbf{x}}$ is not constrained to be $S$-sparse.
The estimation error incurred by an estimator $\hat{\mathbf{x}}(\mathbf{y})$ can be quantified by the mean squared error (MSE)
$\estmsex \triangleq \mathsf{E}_{\mathbf{x}} \big\{ \| \hat{\mathbf{x}}(\mathbf{y}) \rmv-\rmv \mathbf{x} \|^{2}_{2} \big\}$,
where the notation $\mathsf{E}_{\x} \{ \cdot \}$ indicates that the expectation is taken with respect to the pdf $f(\mathbf{y};  \xtrue)$ parameterized by $\x$.
Note that $\estmsex$ 
%% (as well as all other performance measures to be considered later) 
depends on the true parameter value, $\x$. The MSE can be decomposed as
\be
\label{equ_mse_bias_var}
\estmsex \,=\, {\| \estbiasx \|}_2^{2} \ist+\ist \estvarx \,,
\ee
with the estimator bias  
$\estbiasx \ist\triangleq\, \mathsf{E}_{\mathbf{x}} \{ \hat{\mathbf{x}}(\mathbf{y}) \} - \mathbf{x}$ and the estimator variance 
$\estvarx  \ist\triangleq\, \mathsf{E}_{\mathbf{x}} \big\{ \big\| \hat{\mathbf{x}}(\mathbf{y}) \rmv-  \mathsf{E}_{\mathbf{x}} \{ \hat{\mathbf{x}}(\mathbf{y})  \} \big\|^{2} \big\}$.
A standard approach to defining an optimum estimator is to fix the bias, i.e., 
$\estbiasx \ist\stackrel{!}{=}\ist \mathbf{c}(\mathbf{x})$ for all $\mathbf{x} \!\in\! \mathcal{X}$,
%% with some prescribed function $\mathbf{c}(\mathbf{x})$, 
and minimize the variance $v(\hat{\mathbf{x}}(\cdot);\mathbf{x})$ for all $\mathbf{x} \!\in\! \mathcal{X}$ under this bias constraint.
However, in many cases, such a ``uniformly optimum'' estimator does not exist. It is then natural to consider ``locally optimum'' estimators
that minimize $v(\hat{\mathbf{x}}(\cdot);\mathbf{x}_0)$ only at a given parameter value $\x \!=\! \x_0 \!\in\! \mathcal{X}$. 
This approach is taken here. Note that it follows from \eqref{equ_mse_bias_var} that once the bias is fixed, minimizing the variance is equivalent to minimizing the MSE $\varepsilon(\hat{\mathbf{x}}(\cdot); \x_0)$.

The bias constraint $\estbiasx = \mathbf{c}(\mathbf{x})$ can be equivalently written as the mean constraint
\[
\mathsf{E}_{\mathbf{x}} \{ \hat{\mathbf{x}}(\mathbf{y}) \} \ist=\ist \bm{\gamma}(\mathbf{x}) \,, \quad \text{with} \;\, \bm{\gamma}(\mathbf{x}) \ist\triangleq\ist \mathbf{c}(\mathbf{x}) + \mathbf{x} \,.
\]
Thus, we consider the constrained optimization problem 
\be
\label{equ_min_est_var_vector}
\hat{\mathbf{x}}_{\x_0}(\cdot) \eq \arg\min_{\hat{\mathbf{x}}(\cdot) \in \mathcal{B}_{\gamma}} \! v(\hat{\mathbf{x}}(\cdot); \x_0) \,,
\vspace{-3mm}
\ee 
where 
\[
\mathcal{B}_{\bm{\gamma}} \ist\triangleq\ist \big\{ \hat{\mathbf{x}}(\cdot) \ist\big|\ist \mathsf{E}_{\mathbf{x}} \{ \hat{\mathbf{x}}(\mathbf{y}) \} \!=\rmv \bm{\gamma}(\mathbf{x})
\ist ,\ist \forall \, \mathbf{x} \!\in\! \mathcal{X} \big\} \,.
\]
%%  is the set of estimator functions $\hat{\mathbf{x}}(\cdot)$ such that $\estbiasx = \mathbf{c}(\mathbf{x})$ for all $\mathbf{x} \!\in\! \mathcal{X}$.
The minimum variance achieved by the locally optimum estimator $\hat{\mathbf{x}}_{\x_0}(\cdot)$ at $\x_0$ will be denoted as
\[
%% \label{equ_minvar_vector}
V_{\bm{\gamma}}(\mathbf{x}_{0}) \,\triangleq\, v(\hat{\mathbf{x}}_{\x_0}(\cdot); \x_0) 
  \eq \! \min_{\hat{\mathbf{x}}(\cdot) \in \mathcal{B}_{\bm{\gamma}}} \! v(\hat{\mathbf{x}}(\cdot); \x_0) \,.
\]
This is also known as the \emph{Barankin bound} (for the prescribed mean $\bm{\gamma}(\mathbf{x})$) \cite{barankin49}. 
Using RKHS theory, it can be shown that $\hat{\mathbf{x}}_{\x_0}(\cdot)$ exists, i.e., there exists a unique minimum for \eqref{equ_min_est_var_vector}, provided that   %%% BLIBLABLU
%% the set $\mathcal{B}_{\bm{\gamma}}$ contains 
there exists at least one estimator with mean $\bm{\gamma}(\mathbf{x})$ for all $\mathbf{x} \!\in\! \mathcal{X}$ and finite variance at $\mathbf{x}_0$ (see also Section \ref{sec_rkhs}).
For unbiased estimation, i.e., 
%% $\mathbf{c}(\mathbf{x}) \equiv \mathbf{0}$ or 
$\bm{\gamma}(\mathbf{x}) \rmv\equiv\rmv \x$, $\hat{\mathbf{x}}_{\x_0}(\cdot)$ is called a \emph{locally minimum variance unbiased} (LMVU) estimator. 
Unfortunately, $V_{\bm{\gamma}}(\mathbf{x}_{0})$ is difficult to compute in many cases, including the case of the SLM. Lower bounds on $V_{\bm{\gamma}}(\mathbf{x}_{0})$
are, e.g., the CRB and the Hammersley-Chapman-Robbins bound \cite{GormanHero}. 
 
Let $x_k$, $\hat{x}_k(\mathbf{y})$, 
%% $c_k(\mathbf{x})$, 
and $\gamma_k(\mathbf{x})$ denote the $k\ist$th 
entries of $\x$, $\hat{\mathbf{x}}(\mathbf{y})$, 
%% $\mathbf{c}(\mathbf{x})$, 
and $\bm{\gamma}(\mathbf{x})$, respectively.
We have 
%% $\estmsex = \sum_{k=1}^N \estmsexgk$ and 
$\estvarx = \sum_{k=1}^N v(\hat{x}_{k}(\cdot);\mathbf{x})$
with 
%% the component MSE $\estmsexgk \ist\triangleq\, \mathsf{E}_{\mathbf{x}} \big\{ \big[ \hat{x}_k(\mathbf{y}) - x_k \big]^{2} \big\}$
%% and the component variance 
$v(\hat{x}_{k}(\cdot);\mathbf{x}) \ist\triangleq\, \mathsf{E}_{\mathbf{x}} \big\{ \big[ \hat{x}_k(\mathbf{y}) - \mathsf{E}_{\mathbf{x}} \{ \hat{x}_k(\mathbf{y}) \} \big]^{2} \big\}$.
Thus, 
%% the optimization problem 
\eqref{equ_min_est_var_vector} is equivalent to the $N$ scalar optimization problems
\be
\label{equ_min_est_var}
\hat{x}_{\x_0,k}(\cdot) \eq \arg\min_{\hat{x}_{k}(\cdot) \in \mathcal{B}_{\gamma_{k}}} \!\!\rmv v(\hat{x}_{k}(\cdot);\mathbf{x}_0) \,, \quad k = 1,\ldots,N \,,
\vspace{-1.7mm}
\ee 
where
\vspace{-1.5mm}
\[
\mathcal{B}_{\gamma_{k}} \ist\triangleq\ist \big\{ \hat{x}(\cdot) \ist\big|\ist \mathsf{E}_{\mathbf{x}} \{ \hat{x}(\mathbf{y}) \} \!=\rmv \gamma_k(\mathbf{x})\ist ,\ist \forall \, \mathbf{x} \!\in\! \mathcal{X} \big\} \,.
\]
%%  with $\gamma_k(\mathbf{x}) = c_k(\mathbf{x}) + x_k$. 
The minimum variance achieved by $\hat{x}_{\x_0,k}(\cdot)$ at $\x_0$ is denoted \nolinebreak %%%%%%%%  
as
%% \vspace{-1.2mm}
\be
\label{equ_minvar}
V_{\gamma_k}(\mathbf{x}_{0}) \,\triangleq\, v(\hat{x}_{\x_0,k}(\cdot); \x_0) 
  \eq \! \min_{\hat{x}_{k}(\cdot) \in \mathcal{B}_{\gamma_{k}}} \!\!\rmv v(\hat{x}_{k}(\cdot);\mathbf{x}_0) \,.
\vspace{-1.5mm}
\ee
%% Note that 
%% \vspace{-.5mm}
%% $V_{\bm{\gamma}}(\mathbf{x}_{0}) = \sum_{k=1}^N \!V_{\gamma_k}(\mathbf{x}_{0})$.

\subsection{CRB of the Linear Gaussian Model}\label{sec:min_gauss_model}
%%%%%%%%%%%%%%%%%%%%%%%%%%%%%%%%%%%

\vspace{.5mm}

In our further development, we will make use of the CRB for the \emph{linear Gaussian model} (LGM) defined by
\be
\label{equ_linear_gauss_model}
\mathbf{z} \ist=\ist \mathbf{A} \mathbf{s} + \mathbf{n} \,, 
\ee
with the nonrandom parameter $\mathbf{s} \!\in\! 
%% \mathcal{S} \!\subseteq\! 
\mathbb{R}^{S}\!$ (not assumed sparse),
%%  (not assumed to be sparse), 
the\linebreak %%%%%%%%% 
observation $\mathbf{z} \!\in\! \mathbb{R}^{M}\rmv\rmv$, the known matrix $\mathbf{A} \!\in\! \mathbb{R}^{M \times S}\rmv\rmv$, 
and white Gaussian noise $\mathbf{n} \sim \mathcal{N}(\mathbf{0}, \sigma^{2} \mathbf{I})$.
As before, we assume that $M \!\ge\! S$; furthermore, we assume that $\mathbf{A}$ has full column rank, i.e., 
%% the product 
$\mathbf{A}^{\rmv\rmv T} \!\mathbf{A} \!\in\! \mathbb{R}^{S \times S}\rmv$ is nonsingular.
%% Note that the SLM \eqref{equ_observation_model} is a special case because condition \eqref{equ_spark_cond} is satisfied.
The relationship of this model with the SLM, as well as the different notation and different dimension ($S$ instead of $N$), 
will become clear in Section \ref{sec_lower_bounds}.
%% Sections \ref{sec:isometry} and \ref{sec:another_lowerbound}.

Consider estimators $\hat{s}_{k}(\mathbf{z})$ of the $k\ist$th parameter component $s_k$ whose bias is equal to some prescribed differentiable function $\tilde{c}_k(\mathbf{s})$, 
i.e., $b(\hat{s}_{k}(\cdot);\mathbf{s}) = \tilde{c}_k(\mathbf{s})$ or equivalently 
$\mathsf{E}_{\mathbf{s}} \big\{ \hat{s}_{k}(\mathbf{z}) \big\} = \tilde{\gamma}_k(\mathbf{s})$ with $\tilde{\gamma}_k(\mathbf{s}) \triangleq \tilde{c}_{k}(\mathbf{s}) + s_k$,
for all $\mathbf{s} \!\in\! \mathbb{R}^{S}\rmv\rmv$. Let $V_{\tilde{\gamma}_{k}}^{\text{LGM}}(\mathbf{s}_{0})$ denote the minimum variance achievable by such estimators at a given true parameter $\mathbf{s}_{0}$.
The CRB $C^{\text{LGM}}_{\tilde{\gamma}_k} (\mathbf{s}_{0})$ is the following lower bound on the minimum variance \cite{kay}: 
\be
\label{equ_CRB_LGM}
%% v(\hat{s}_{k}(\cdot);\mathbf{s}) 
V_{\tilde{\gamma}_{k}}^{\text{LGM}}(\mathbf{s}_{0}) \ist\geq\, %% {\mbox{CRB}}
C^{\text{LGM}}_{\tilde{\gamma}_k} (\mathbf{s}_{0})
\,\triangleq\, \sigma^{2} \,\tilde{\mathbf{r}}_k^{T}\! (\mathbf{s}_0) \ist {( \mathbf{A}^{\rmv\rmv T} \!\mathbf{A} )}^{\!-1} \rmv\tilde{\mathbf{r}}_k(\mathbf{s}_{0}) \,,
\ee
where 
%% $\mathbf{r}_k(\mathbf{s}) \triangleq \frac{\partial \ist \tilde{\gamma}_k(\mathbf{s}) }{\partial \mathbf{s}}$,
$\tilde{\mathbf{r}}_k(\mathbf{s}) \triangleq \partial \ist \tilde{\gamma}_k(\mathbf{s}) /\partial \mathbf{s}$,
i.e., $\tilde{\mathbf{r}}_k(\mathbf{s})$ is the vector of dimension $S$ whose $l\ist$th entry is $\partial \ist \tilde{\gamma}_k(\mathbf{s}) / \partial s_l$. 
%% = \partial \ist (c_{k}(\mathbf{s}) \rmv+\rmv s_k ) / \partial s_l 
%% = \partial \ist \tilde{c}_{k}(\mathbf{s}) / \partial s_l + \delta_{k,l}$. 
We note that $V_{\tilde{\gamma}_{k}}^{\text{LGM}}(\mathbf{s}_{0})  = 
C^{\text{LGM}}_{\tilde{\gamma}_k} (\mathbf{s}_{0})$ if $\tilde{\gamma}_k(\mathbf{s})$ is an affine function of $\mathbf{s}$.
In particular, this includes the unbiased case ($\tilde{\gamma}_k(\mathbf{s}) \rmv\equiv\rmv s_k$).

\vspace{-.5mm}

%% \newpage %%%%%%%%%%

%%%%%%%%%%%%%%%%%%%%%%%%%%%%%%%%%%%
\section{The RKHS Framework}\label{sec_rkhs}
%%%%%%%%%%%%%%%%%%%%%%%%%%%%%%%%%%%

\vspace{.8mm}

In this section, we review some RKHS fundamentals which will provide a basis for our further development. 
Consider a set 
%% \emph{a priori}, 
$\mathcal{X}$ %%% BLIBLABLU
(not necessarily a linear
space) and a positive 
semidefinite\footnote{That %%%%%%%%%
%% A function $R(\mathbf{x}, \mathbf{x}') \!: \mathcal{X} \!\times\! \mathcal{X} \rightarrow \mathbb{R}$ is called positive semidefinite if 
is, for any finite set $\{ \x_{k} \}_{k=1,\ldots,P}$ with $\x_{k} \!\in\! \mathcal{X}$, the matrix $\mathbf{R} \!\in\! \mathbb{R}^{P \times P}$ with entries ${( \mathbf{R} )}_{k,l} \rmv\rmv\triangleq R(\x_{k}, \x_{l})$ is positive 
semidefinite.} %%%%%%%%%%%
``kernel'' function $R(\mathbf{x}, \mathbf{x}') \!: \mathcal{X} \!\times\! \mathcal{X} \rmv\rmv\rightarrow \mathbb{R}$. 
For each fixed $\mathbf{x}' \!\rmv \in\! \mathcal{X}$, the function $f_{\mathbf{x}'}(\x) \triangleq R(\x, \mathbf{x}')$ maps $\mathcal{X}$ into $\mathbb{R}$.
The RKHS $\mathcal{H}(R)$ is a Hilbert space of functions $f \!: \mathcal{X} \!\rightarrow\rmv \mathbb{R}$ which is defined as the 
closure of the linear span of the set of functions ${ \{ f_{\mathbf{x}'}(\x) = R(\x, \mathbf{x}') \} }_{\mathbf{x}' \in \mathcal{X}}$.
This closure is taken with respect to the topology given by the scalar product ${\langle \cdot\ist\ist , \cdot \rangle}_{\mathcal{H}(R)}$ which is defined via the \emph{reproducing property} \cite{aronszajn1950} %%% BLIBLABLU
\[
\big\langle f(\cdot) , R(\cdot,\mathbf{x}') \big\rangle_{\mathcal{H}(R)} \,=\, f(\mathbf{x}') \,.
\]
This relation holds for all $f \in \mathcal{H}(R)$ and $\mathbf{x}' \!\rmv \in\! \mathcal{X}$. 
The associated norm is given by 
\pagebreak %%%%%%%%%%
${\| f \|}_{\mathcal{H}(R)} = {\langle f , f \rangle}_{\mathcal{H}(R)}^{1/2}$.

We now consider the constrained 
optimization problem \eqref{equ_min_est_var} for a given mean function $\gamma(\mathbf{x})$ (formerly denoted by $\gamma_k(\mathbf{x})$; we temporarily
drop the subscript $k$ for better readability). According to \cite{Parzen59,Duttweiler73b}, for certain classes of parametrized pdf's 
$f(\mathbf{y}; \mathbf{x})$ (which include the Gaussian pdf in \eqref{equ_pdf}), 
one can associate with this optimization problem an RKHS $\mathcal{H}(R_{\mathbf{x}_{0}})$ 
whose kernel $R_{\mathbf{x}_{0}}(\mathbf{x}, \mathbf{x}') \!: \mathcal{X} \!\times\! \mathcal{X} \!\rightarrow\! \mathbb{R}$ is given by 
\begin{align*} 
R_{\mathbf{x}_{0}}(\mathbf{x}, \mathbf{x}' ) &\,\triangleq\, \mathsf{E}_{\mathbf{x}_{0}} \rmv\bigg \{ \frac{f(\mathbf{y}; \mathbf{x})}{f(\mathbf{y}; \mathbf{x}_{0} )} \frac{f(\mathbf{y}; \mathbf{x}')}{f(\mathbf{y};\mathbf{x}_{0})} \bigg \}  \\[1.5mm] 
&\eq\! \int_{\mathbb{R}^M} \!\! \frac{ f(\mathbf{y}; \mathbf{x}) \, f(\mathbf{y}; \mathbf{x}') }{ f(\mathbf{y}; \mathbf{x}_{0}) } \, d\mathbf{y} \,.\\[-4mm] 
&
\end{align*} 
%% where $\mathbf{x}_{0} \in \mathcal{X}$. 
It can be shown \cite{Parzen59,Duttweiler73b} that $\gamma(\mathbf{x}) \!\in\! \mathcal{H}(R_{\mathbf{x}_{0}})$
%% ---i.e., the set $\mathcal{B}_{\gamma}$ is nonempty---
if and only if there exists at least one estimator with mean $\gamma(\mathbf{x})$ for all $\mathbf{x}$ and finite variance at $\mathbf{x}_{0}$.  %%% BLIBLABLU
Furthermore, under this condition, the minimum variance 
%% achieved by the optimum estimator, 
$V_{\gamma}(\mathbf{x}_{0})$ in \eqref{equ_minvar} is finite and
%% ---which is finite---
allows the following expression involving the norm 
\vspace{-.5mm}
${\| \gamma \|}_{\mathcal{H}(R_{\mathbf{x}_{0}})}$:
\begin{equation} 
V_{\gamma}(\mathbf{x}_{0}) \eq {\| \gamma \|}^{2}_{\mathcal{H}(R_{\mathbf{x}_{0}})} \rmv-\ist \gamma^2(\mathbf{x}_{0}) \,. 
\label{equ_min_expr_RKHS}
\end{equation}
This is an RKHS formulation of the Barankin bound. Unfortunately, the norm ${\| \gamma \|}_{\mathcal{H}(R_{\mathbf{x}_{0}})}$
is often difficult to compute.

For the SLM in \eqref{equ_observation_model}, \eqref{equ_SLM_parameter}, \eqref{equ_pdf},
%% $R_{\mathbf{x}_{0}}(\mathbf{x},\mathbf{x}')$ with $\mathbf{x}_{0} \!\in\! \mathcal{X}_{S}$
%% we have 
$\mathcal{X} \!=\! \mathcal{X}_{S}$; the kernel here is a mapping $\mathcal{X}_{S} \rmv\times\rmv \mathcal{X}_{S} \rightarrow \mathbb{R}$ which is easily shown to be given 
%% \vspace{-1mm}
by 
\be
\label{equ_kernel_SLM}
R_{\mathbf{x}_{0}}(\mathbf{x},\mathbf{x}') 
  \eq \exp \rmv\rmv \bigg( \frac{1}{\sigma^{2}} \ist ( \mathbf{x} \!-\! \mathbf{x}_{0} )^{T} \ist\mathbf{H}^{T} \mathbf{H} \ist\ist ( \mathbf{x}' \!\!-\! \mathbf{x}_{0} ) \bigg) \,,
\ee
where $\mathbf{x}_{0} \!\in\! \mathcal{X}_{S}$. An RKHS can also be defined for the LGM in \eqref{equ_linear_gauss_model}. Here, $\mathcal{X} \rmv=\rmv \mathbb{R}^{S}$, 
and the kernel $R_{\mathbf{s}_{0}}^{\text{LGM}}(\mathbf{s},\mathbf{s}')$ with $\mathbf{s}_{0} \!\in\! \mathbb{R}^{S}$ is 
a mapping $\mathbb{R}^{S} \!\rmv\times\rmv \mathbb{R}^{S} \!\rightarrow \mathbb{R}$ given by 
\be
\label{equ_kernel_LGM}
R_{\mathbf{s}_{0}}^{\text{LGM}}(\mathbf{s},\mathbf{s}') 
  \eq \exp \rmv\rmv \bigg( \frac{1}{\sigma^{2}} \ist ( \mathbf{s} \!-\! \mathbf{s}_{0} )^{T} \rmv\rmv\mathbf{A}^{\rmv\rmv T} \!\mathbf{A} \ist ( \mathbf{s}' \!\!-\! \mathbf{s}_{0} ) \bigg) \,.
\ee
%% The corresponding RKHS will be denoted by $\mathcal{H}(R_{\mathbf{s}_{0}}^{\text{LGM}})$.
%% The kernels $R_{\mathbf{x}_{0}}(\mathbf{x},\mathbf{x}')$ and $R_{\mathbf{s}_{0}}^{\text{LGM}}(\mathbf{s},\mathbf{s}')$ 
%% are seen to have the same functional form. However, they 
Note that these kernels differ in their domain, which is $\mathcal{X}_{S} \rmv\times\rmv \mathcal{X}_{S}$ for $R_{\mathbf{x}_{0}}(\mathbf{x},\mathbf{x}')$ 
and $\mathbb{R}^{S} \!\rmv\times\rmv \mathbb{R}^{S}$ for $R_{\mathbf{s}_{0}}^{\text{LGM}}(\mathbf{s},\mathbf{s}')$.

%% \newpage %%%%%%%%%

\vspace{-.2mm}

%%%%%%%%%%%%%%%%%%%%%%%%%%%%%%%%%%%
\section{A Lower Bound on the Estimator Variance} 
\label{sec_lower_bounds}
%%%%%%%%%%%%%%%%%%%%%%%%%%%%%%%%%%%

\vspace{.8mm}

%% \subsection{Lower Bounds based on the RKHS framework}

We now continue our treatment of the SLM estimation problem.
In what follows, $V_{\gamma}(\mathbf{x}_{0})$ will be understood to denote the bias-constrained minimum 
%% component 
variance 
%% $V_{\gamma}(\mathbf{x}_{0})$ (see 
\eqref{equ_minvar} \emph{specifically for the SLM}. This means, in particular, that $\mathcal{X} \!=\! \mathcal{X}_{S}$, and hence the set of admissible
estimators is given by
\be
\label{equ_est-set}
\mathcal{B}_{\gamma} =\ist \big\{ \hat{x}(\cdot) \,\big|\, \mathsf{E}_{\mathbf{x}} \{ \hat{x}(\mathbf{y}) \} \rmv=\rmv \gamma(\mathbf{x}) \ist , \ist \forall \, \mathbf{x} \!\in\! \mathcal{X}_{S} \big\} \,.
\ee 
We will next
%% use the RKHS framework to 
%% (with $\mathcal{X} \!=\! \mathcal{X}_{S}$) 
derive a lower bound on $V_{\gamma}(\mathbf{x}_{0})$.
%% denoted
%% \be
%% \label{equ_min_var_SLM} 
%% V^{\text{SLM}}_{\gamma}(\mathbf{x}_{0}) \triangleq \min_{\hat{x}(\cdot) \in\mathcal{B}_{\gamma}}  \!\rmv v(\hat{x}(\cdot);\mathbf{x}_0),
%% \ee
%% for the SLM, i.e., the parameter set $\mathcal{X}$ in \eqref{equ_min_est_var_recall} is given by $\mathcal{X}_{S}$.
%% where, as before, $\mathcal{B} = \big\{ \hat{x}(\cdot) \,\big|\, \mathsf{E}_{\mathbf{x}} \{ \hat{x}(\mathbf{y}) \} 
%% = \gamma(\mathbf{x}) \ist , \,\ist \forall \, \mathbf{x} \!\in\! \mathcal{X}_{S} \big\}$ 
%% with $\gamma(\mathbf{x}) = c(\mathbf{x}) + x$. 

\vspace{-1.5mm}

%% \subsection{Lower Bound on \emph{$V_{\gamma}(\mathbf{x}_{0})$:} Relaxing the Bias Constraint}\label{sec:relax}
\subsection{Relaxing the Bias Constraint}\label{sec:relax}
%%%%%%%%%%%%%%%%%%%%%%%%%%%%%%%%%%%
%%  with Prescribed Bias

\vspace{.5mm}

The first step in this derivation is to relax the bias constraint $\hat{x}(\cdot) \!\in\! \mathcal{B}_{\gamma} $. 
Let $\mathcal{K} \triangleq \{ k_{1}, 
%% k_{2}, ..., 
\ldots, k_{S} \}$ be a fixed set of $S$ different indices $k_i \in \{ 1, \ldots, N \}$ (not related to $\rm{supp}(\x_0)$),
%% . Furthermore let $\mathcal{X}_{S}^{\mathcal{K}}$ denote 
%% the set of all $\mathbf{x} \!\in\! \mathcal{X}_{S}$ whose support (i.e., the set of indices of the nonzero components) equals $\mathcal{K}$:
and let
\[
%% \label{equ_def_subset_U}
\mathcal{X}_{S}^{\mathcal{K}} \ist\triangleq\ist \{ \mathbf{x} \!\in\! \mathcal{X}_{S} \ist | \ist \supp(\mathbf{x}) \rmv\subseteq\rmv \mathcal{K} \} \,.
\]
Clearly, $\mathcal{X}_{S}^{\mathcal{K}} \!\subseteq\! \mathcal{X}_{S}$; however, contrary to $\mathcal{X}_{S}$, $\mathcal{X}_{S}^{\mathcal{K}}$ 
%% has the important property of being 
is a linear subspace of $\mathbb{R}^{N}\rmv\rmv$. 
Let $\mathcal{B}_{\gamma}^{\mathcal{K}}$ be the set of all 
\pagebreak %%%%%%%%%
estimators with mean $\gamma(\mathbf{x})$ for all $\x \!\in\! \mathcal{X}_{S}^{\mathcal{K}}$
(but not necessarily for all $\x \!\in\! \mathcal{X}_{S}$), 
%% \vspace*{-4.5mm}
i.e.,
\[
\mathcal{B}_{\gamma}^{\mathcal{K}} \triangleq\ist \big\{ \hat{x}(\cdot) \,\big|\, \mathsf{E}_{\mathbf{x}} \{ \hat{x}(\mathbf{y}) \} \rmv=\rmv \gamma(\mathbf{x}) \ist , \ist \forall \, \mathbf{x} \!\in\! \mathcal{X}_{S}^{\mathcal{K}} \big\} \,.
\vspace{-.3mm}
\] 
Comparing with \eqref{equ_est-set}, we see that $\mathcal{B}_{\gamma}^{\mathcal{K}} \!\supseteq\! \mathcal{B}_{\gamma}$. 

Let us now consider the minimum variance among all estimators in $\mathcal{B}_{\gamma}^{\mathcal{K}}$, i.e.,
\be
\label{equ_Vmin_K}
V_{\gamma}^{\mathcal{K}}(\mathbf{x}_{0}) \,\triangleq\rmv \min_{\hat{x}(\cdot) \in \mathcal{B}_{\gamma} ^{\mathcal{K}}} \rmv v(\hat{x}(\cdot);\mathbf{x}_0) \,.
\vspace{-1mm}
\ee
Because $\hat{x}(\cdot) \!\in\! \mathcal{B}_{\gamma}^{\mathcal{K}}$ is a less restrictive constraint than 
%% the constraint 
$\hat{x}(\cdot) \!\in\! \mathcal{B}_{\gamma}$
used in the definition of $V_{\gamma}(\mathbf{x}_{0})$, we have
\be
\label{equ_V_inequ}
V_{\gamma}(\mathbf{x}_{0})  \ist\ge V_{\gamma}^{\mathcal{K}}(\mathbf{x}_{0}) \,,
\vspace{-.3mm}
\ee
i.e., $V_{\gamma}^{\mathcal{K}}(\mathbf{x}_{0})$ is a lower bound on $V_{\gamma}(\mathbf{x}_{0})$. 
A closed-form expression of $V_{\gamma}^{\mathcal{K}}(\mathbf{x}_{0})$ appears to be difficult to obtain in the general case, because 
$\mathbf{x}_{0} \!\not\in\! \mathcal{X}_{S}^{\mathcal{K}}$ in general.
Therefore, we will use RKHS theory to derive a lower bound on $V_{\gamma}^{\mathcal{K}}(\mathbf{x}_{0})$. 
%% This will be done by using RKHS theory and establishing a relation with the CRB for the LGM in \eqref{equ_CRB_LGM}.

\vspace{-1.8mm}

%% \newpage %%%%%%%%

\subsection{Two Isometric RKHSs}\label{sec:isometry}
%%%%%%%%%%%%%%%%%%%%%%%%%%%%%%%%%%%
%%  with Prescribed Bias

\vspace{.7mm}

An RKHS 
%% $\RKHSc$ 
for the SLM can also be defined on $\mathcal{X}_{S}^{\mathcal{K}}$, using a kernel 
%% $R_{\mathbf{x}_{0}}^{\mathcal{K}}(\mathbf{x},\mathbf{x}')$ 
$R_{\mathbf{x}_{0}}^{\mathcal{K}} \!: \mathcal{X}_{S}^{\mathcal{K}} \rmv\times\rmv \mathcal{X}_{S}^{\mathcal{K}} \rightarrow \mathbb{R}$ that is given by 
the right-hand side of \eqref{equ_kernel_SLM} but whose arguments $\mathbf{x},\mathbf{x}'$ are assumed to be in $\mathcal{X}_{S}^{\mathcal{K}}$ 
and not just in $\mathcal{X}_{S}$ (however, recall that $\mathbf{x}_{0} \!\not\in\!\rmv \mathcal{X}_{S}^{\mathcal{K}}\rmv\rmv$ in general).
%% \be
%% \label{equ_kernel_K}
%% R_{\mathbf{x}_{0}}^{\mathcal{K}}(\mathbf{x},\mathbf{x}') 
%%   \eq \exp \rmv\rmv \bigg( \frac{1}{\sigma^{2}} \ist ( \mathbf{x} \!-\! \mathbf{x}_{0} )^{T} \ist\mathbf{H}^{T} \mathbf{H} \ist\ist ( \mathbf{x}' \!\!-\! \mathbf{x}_{0} ) \bigg) \,.
%% \ee
%% Note that whereas $\mathbf{x},\mathbf{x}' \in \mathcal{X}_{S}^{\mathcal{K}}$, the true parameter $\mathbf{x}_{0}$ 
%% is not an element of $\mathcal{X}_{S}^{\mathcal{K}}$ in general. 
This RKHS will be denoted $\RKHSc$. The minimum variance $V_{\gamma}^{\mathcal{K}}(\mathbf{x}_{0})$ 
in \eqref{equ_Vmin_K} can then be expressed as (cf.\ \eqref{equ_min_expr_RKHS})
\be
\label{equ_min_expr_RKHS_alt}
V_{\gamma}^{\mathcal{K}}(\mathbf{x}_{0}) \eq {\| \gamma \|}^{2}_{\RKHSc } \rmv-\ist \gamma^{2}(\mathbf{x}_{0}) \,.
\ee

In order to develop this expression, we define some notation. 
Consider an 
%% ordered 
index set $\mathcal{I} = \{k_1,\ldots,k_{|\mathcal{I}|}\} \subseteq \{1,\ldots,N\}$. 
We denote by $\mathbf{H}_{\mathcal{I}} \!\in\! \mathbb{R}^{M \times |\mathcal{I}|}\rmv$ the 
%% $M \!\times\! |\mathcal{I}|$ 
submatrix of our 
%% $M \!\times\! N$ 
matrix $\mathbf{H} \!\in\! \mathbb{R}^{M \times N}\!$
whose $i\ist$th column is given by the $k_{i}\ist$th column of $\mathbf{H}$. Furthermore, for a vector $\mathbf{x} \!\in\! \mathbb{R}^{N}\rmv\rmv$,
we denote by $\mathbf{x}^{\mathcal{I}} \!\!\in\! \mathbb{R}^{|\mathcal{I}|}\rmv$ the subvector whose $i\ist$th entry is the $k_{i}\ist$th entry of $\mathbf{x}$. 

We now introduce a second RKHS. Consider the LGM in \eqref{equ_linear_gauss_model} with 
%% system 
matrix $\mathbf{A} \rmv= \mathbf{H}_{\mathcal{K}} \in \mathbb{R}^{M \times S}\rmv\rmv$, 
and let $\mathcal{H}(R_{\mathbf{s}_{0}}^{\text{LGM}})$ with $\mathbf{s}_{0} \!\in\! \mathbb{R}^S\rmv$ denote the RKHS for that LGM as defined by the kernel $R_{\mathbf{s}_{0}}^{\text{LGM}} \!: \mathbb{R}^{S} \!\times\rmv \mathbb{R}^{S} \rmv\rmv\rightarrow \mathbb{R}$ in \eqref{equ_kernel_LGM}.
Exploiting the linear-subspace structure of $\mathcal{X}_{S}^{\mathcal{K}}$, it can be shown that our RKHS $\RKHSc$ for a given $\x_0$
%% defined by \eqref{equ_kernel_K} 
is \emph{isometric} to $\mathcal{H}(R_{\mathbf{s}_{0}}^{\text{LGM}})$ with $\mathbf{s}_{0}$ chosen as
%% , i.e., the RKHS associated to the optimization problem \eqref{equ_min_est_var} for a LGM 
\begin{equation} 
\mathbf{s}_{0} \rmv\eq \mathbf{H}_{\mathcal{K}}^\dag \ist \mathbf{H} \, \mathbf{x}_{0} \, .
  %%  \in \mathbb{R}^S . 
\label{equ_s0_expr}
\end{equation}
Here, $\mathbf{H}_{\mathcal{K}}^\dag \triangleq ( \mathbf{H}_{\mathcal{K}}^{T} \mathbf{H}_{\mathcal{K}} )^{-1} \mathbf{H}_{\mathcal{K}}^{T} 
\in \mathbb{R}^{S \times M}\!$ is the pseudo-inverse of $\mathbf{H}_{\mathcal{K}}$ (recall that $M \!\ge\! S$, and note that 
$( \mathbf{H}_{\mathcal{K}}^{T} \mathbf{H}_{\mathcal{K}} )^{-1}$ is guaranteed to exist because of our assumption \eqref{equ_spark_cond}).
More specifically, the isometry $\mathsf{J} \!: \RKHSc \rightarrow \mathcal{H}(R_{\mathbf{s}_{0}}^{\text{LGM}})$ mapping each $f \!\in\rmv \RKHSc$ to an $\tilde{f} \rmv\in\rmv \mathcal{H}(R_{\mathbf{s}_{0}}^{\text{LGM}})$
%% , $J: \RKHSc \rightarrow  \mathcal{H}(R_{\mathbf{s}_{0}}^{\text{LGM}}) : f \mapsto f'$ 
is given by 
\begin{equation} 
%f'(h^{-1}_{\mathcal{K}}(\mathbf{x}) )   = f(\mathbf{x}) z 
\mathsf{J} \{ f(\mathbf{x}) \} \ist=\ist \tilde{f}(\mathbf{x}^{\mathcal{K}}) \ist=\ist \beta_{\mathbf{x}_{0}} \ist f(\mathbf{x}) \,, \quad\; 
%% \text{for all} \;\, 
\mathbf{x} \!\in\! \mathcal{X}_{S}^{\mathcal{K}} \,,
\vspace{-1.5mm}
\label{equ_J_expr}
\end{equation}
where 
\be
\label{equ_beta-def}
\beta_{\mathbf{x}_{0}} \rmv\triangleq\, \exp \rmv\rmv \bigg( \!\!\rmv-\rmv \frac{1}{2 \sigma^{2}} \rmv \big\| (\mathbf{I} \!-\! \mathbf{P}_{\!\mathcal{K}})
  \, \mathbf{H} \mathbf{x}_{0}  \big\|_2^{2} \bigg) \ist\ist .
%% \label{equ_beta_def}
\vspace{-1mm}
\ee
Here, $\mathbf{P}_{\!\mathcal{K}} \triangleq \mathbf{H}_{\mathcal{K}} \mathbf{H}_{\mathcal{K}}^\dag$ is the orthogonal projection matrix on the range of $\mathbf{H}_{\mathcal{K}}$. 
The factor $\beta_{\mathbf{x}_{0}}\!$ can be interpreted as a 
measure of the distance between the point $\mathbf{H} \mathbf{x}_{0}$ and the subpsace $\mathcal{X}_{S}^{\mathcal{K}}$ 
\pagebreak %%%%%%%%%
associated with the index set $\mathcal{K}$. We can write \eqref{equ_J_expr} as
\[
%% \mathsf{J} \{ f(\mathbf{x}) \} \eq 
\tilde{f}(\mathbf{s}) \ist=\ist \beta_{\mathbf{x}_{0}} \ist f(\mathbf{x}(\mathbf{s})) \,, \quad\; 
%% \text{for all} \;\, 
\mathbf{s} \!\in\! \mathbb{R}^{S} ,
\]
where $\mathbf{x}(\mathbf{s})$ denotes the $\mathbf{x} \!\in\! \mathcal{X}_{S}^{\mathcal{K}}$ for which $\mathbf{x}^{\mathcal{K}} \!=\rmv \mathbf{s}$
(i.e., the $S$ entries of $\mathbf{s}$ appear in $\mathbf{x}(\mathbf{s})$ at the appropriate positions within $\mathcal{K}$,
and the $N\!-\!S$ remaining entries of $\mathbf{x}(\mathbf{s})$ are zero). 
%% (i.e., the $S$ entries of $\mathbf{x}(\mathbf{s})$ located within the index set $\mathcal{K}$ are equal to the appropriate entries of $\mathbf{s}$
%% and the $N\!-\!S$ remaining entries of $\mathbf{x}(\mathbf{s})$ are zero). 

Consider now the image of 
%% the prescribed mean function 
$\gamma(\mathbf{x})$ under the mapping $\mathsf{J}$,
\be
\label{equ_relation_tilde_gamma}
\tilde{\gamma}(\mathbf{s}) \ist\triangleq\ist \mathsf{J}\{ \gamma(\mathbf{x}) \} \ist=\ist \beta_{\mathbf{x}_{0}} \gamma(\mathbf{x}(\mathbf{s})) \,, \quad\; 
%% \text{for all} \;\, 
\mathbf{s} \!\in\! \mathbb{R}^{S} .
\ee
Since $\mathsf{J}$ is an isometry, we have ${\| \tilde{\gamma} \|}^{2}_{\mathcal{H}(R_{\mathbf{s}_{0}}^{\text{LGM}})} \rmv\rmv\!=\! {\|\gamma \|}^{2}_{\RKHSc}$.
Combining this identity with \eqref{equ_min_expr_RKHS_alt}, we obtain
\be
\label{equ_V_gamma-tilde}
V_{\gamma}^{\mathcal{K}}(\mathbf{x}_{0}) \eq {\| \tilde{\gamma} \|}^{2}_{\mathcal{H}(R_{\mathbf{s}_{0}}^{\text{LGM}})} - \gamma^{2}(\mathbf{x}_{0}) \,.
\vspace{-2mm}
\ee

%% \newpage %%%%%%%%

\subsection{Lower Bound on \emph{$V_{\gamma}^{\mathcal{K}}(\mathbf{x}_{0})$}} \label{sec:another_lowerbound}
%%%%%%%%%%%%%%%%%%%%%%%%%%%%%%%%%%%
%%  with Prescribed Bias

\vspace{.5mm}

We will now use expression \eqref{equ_V_gamma-tilde} to derive a lower bound on $V_{\gamma}^{\mathcal{K}}(\mathbf{x}_{0})$ in terms of the CRB for the LGM in \eqref{equ_CRB_LGM}.
Consider the minimum estimator variance for the LGM under the constraint of the prescribed mean function $\tilde{\gamma}(\mathbf{s})$, 
$V_{\tilde{\gamma}}^{\text{LGM}}(\mathbf{s}_{0})$,
still for $\mathbf{A} \rmv= \mathbf{H}_{\mathcal{K}}$ and for $\mathbf{s}_{0}$ given by \eqref{equ_s0_expr}.
We have (cf.\ \eqref{equ_min_expr_RKHS})
\[
%% \label{equ_V_gama-tilde_LGM}
V_{\tilde{\gamma}}^{\text{LGM}}(\mathbf{s}_{0}) \eq {\| \tilde{\gamma} \|}^{2}_{\mathcal{H}(R_{\mathbf{s}_{0}}^{\text{LGM}})} - \tilde{\gamma}^{2}(\mathbf{s}_{0}) \,.
\]
Combining with \eqref{equ_V_gamma-tilde}, we obtain the relation
\[
%% \label{equ_V_V}
V_{\gamma}^{\mathcal{K}}(\mathbf{x}_{0}) \eq V_{\tilde{\gamma}}^{\text{LGM}}(\mathbf{s}_{0}) \ist+\ist \tilde{\gamma}^{2}(\mathbf{s}_{0}) - \gamma^{2}(\mathbf{x}_{0}) \,.
\]
Using the CRB $V_{\tilde{\gamma}}^{\text{LGM}}(\mathbf{s}_{0}) \geq C^{\text{LGM}}_{\tilde{\gamma}} (\mathbf{s}_{0})$ (see \eqref{equ_CRB_LGM}) yields 
%% the lower bound
\be
\label{equ_bound_shorthand}
V_{\gamma}^{\mathcal{K}}(\mathbf{x}_{0}) \ist\ge L_{\gamma}^{\mathcal{K}}\rmv (\mathbf{x}_{0}) \,,
\vspace{-2mm}
\ee
with
\vspace{1mm}
\be
\label{equ_V_CRB_Ldef}
L_{\gamma}^{\mathcal{K}}\rmv (\mathbf{x}_{0}) \,\triangleq\, C^{\text{LGM}}_{\tilde{\gamma}} (\mathbf{s}_{0}) \ist+\ist \tilde{\gamma}^{2}(\mathbf{s}_{0}) - \gamma^{2}(\mathbf{x}_{0}) \,.
\ee
Finally, using \eqref{equ_relation_tilde_gamma}
%% , i.e., $\tilde{\gamma}(\mathbf{s}) = \beta_{\mathbf{x}_{0}} \gamma(\mathbf{x}(\mathbf{s}))$, 
and the implied CRB relation $C^{\text{LGM}}_{\tilde{\gamma}} (\mathbf{s}_{0}) 
= \beta_{\mathbf{x}_{0}}^2 C^{\text{LGM}}_{\gamma(\mathbf{x}(\mathbf{s}))}(\mathbf{s}_{0})$, the lower bound \eqref{equ_V_CRB_Ldef} can be reformulated as
\be
\label{equ_V_CRB_Lfinal}
L_{\gamma}^{\mathcal{K}}\rmv (\mathbf{x}_{0}) \eq \beta_{\mathbf{x}_{0}}^2 \big[ C^{\text{LGM}}_{\gamma(\mathbf{x}(\mathbf{s}))}(\mathbf{s}_{0}) 
  \ist+\ist \gamma^2(\mathbf{x}(\mathbf{s}_0)) \big] - \gamma^{2}(\mathbf{x}_{0}) \,.
\ee
Here, $C^{\text{LGM}}_{\gamma(\mathbf{x}(\mathbf{s}))}(\mathbf{s}_{0})$ 
denotes the CRB for prescribed mean function $\gamma'(\mathbf{s}) = \gamma(\mathbf{x}(\mathbf{s}))$, which is given by (see \eqref{equ_CRB_LGM})
\be
\label{equ_CRB_LGM_H}
C^{\text{LGM}}_{\gamma(\mathbf{x}(\mathbf{s}))}(\mathbf{s}_{0})
\eq \sigma^{2} \,\mathbf{r}^{T}\! (\mathbf{s}_0) \ist {( \mathbf{H}_{\mathcal{K}}^{T} \mathbf{H}_{\mathcal{K}} )}^{\!-1} \rmv\mathbf{r}(\mathbf{s}_{0}) \,,
\ee
where $\mathbf{r}(\mathbf{s}) \triangleq \partial \ist \gamma(\mathbf{x}(\mathbf{s})) / \partial \mathbf{s}$ and $\mathbf{s}_{0}$ is related to $\mathbf{x}_{0}$ via \eqref{equ_s0_expr}.

To summarize, we have the following chain of lower bounds on the bias-constrained variance at $\mathbf{x}_{0}$:
%%  of any estimator $\hat{x}(\cdot)$ whose mean equals $\gamma(\mathbf{x})$ for all 
%%  \vspace*{-1mm}
%%  $\mathbf{x} \!\in\! \mathcal{X}_{S}$:
\be
\label{equ_general_lower_bound_g_k}
v(\hat{x}(\cdot);\mathbf{x}_{0}) 
\ist\stackrel{\eqref{equ_minvar}}{\geq}\ist V_{\gamma}(\mathbf{x}_{0}) 
\ist\stackrel{\eqref{equ_V_inequ}}{\geq}\ist V_{\gamma}^{\mathcal{K}}(\mathbf{x}_{0}) 
\ist\stackrel{\eqref{equ_bound_shorthand}}{\geq}\ist L_{\gamma}^{\mathcal{K}}\rmv (\mathbf{x}_{0}) \,.
\ee
While $L_{\gamma}^{\mathcal{K}}\rmv (\mathbf{x}_{0})$ is the loosest of these bounds, it is attractive because
of its closed-form expression in \eqref{equ_V_CRB_Lfinal} (together with \eqref{equ_CRB_LGM_H} and \eqref{equ_s0_expr}).
We note that the inequality \eqref{equ_bound_shorthand} becomes an equality if $\tilde{\gamma}(\mathbf{s})$ is an affine function of $\mathbf{s}$,
or equivalently (see \eqref{equ_relation_tilde_gamma}), if $\gamma(\mathbf{x})$ is an affine function of $\mathbf{x}$. 
In particular, this includes the unbiased case 
\pagebreak %%%%%%%%%
($\gamma(\mathbf{x}) \rmv\equiv\rmv x$).

Recalling that 
%% the variance of the vector estimator $\hat{\mathbf{x}}(\cdot)$ is given by
$v(\hat{\mathbf{x}}(\cdot);\mathbf{x}_{0}) = \sum_{k=1}^N v(\hat{x}_{k}(\cdot);\mathbf{x}_0)$ (we now reintroduce the subscript $k$),
a lower bound on $v(\hat{\x}(\cdot);\mathbf{x}_{0})$ is obtained from \eqref{equ_general_lower_bound_g_k} 
\vspace{-1.7mm}
as
\[
%% \label{equ_lower_vec_est_nontightened}
 v(\hat{\mathbf{x}}(\cdot);\mathbf{x}_{0}) \ist\geq\ist \sum_{k=1}^{N} L_{\gamma_k}^{\mathcal{K}_k}\rmv (\mathbf{x}_{0}) \,.
\vspace{-.7mm}
\] 
For a high lower bound, the index sets $\mathcal{K}_k$ should in general be chosen such that the respective factors $\beta_{\mathbf{x}_{0},k}^2$ in \eqref{equ_V_CRB_Lfinal}
are large. (This means that the ``distances'' between $\mathbf{H} \mathbf{x}_{0}$ and $\mathcal{X}_{S}^{\mathcal{K}_k}$ are small, see \eqref{equ_beta-def}.)
%% Finally, since the index sets $\mathcal{K}_k$ have not been specified so far, we can tighten 
%% this bound by maximizing each $L_{\gamma_k}^{\mathcal{K}_k}\rmv (\mathbf{x}_{0})$ with respect to $\mathcal{K}_k$.
%% , i.e., by calculating $L_{\gamma_k}^*\rmv (\mathbf{x}_{0}) \triangleq \max_{\mathcal{K}_k} L_{\gamma_k}^{\mathcal{K}_k}\rmv (\mathbf{x}_{0})$. 
Formally using the optimum $\mathcal{K}_k$ for each $k$, we arrive at the main result of this 
\vspace{2mm}
paper.
 
\noindent 
%% \begin{theorem} 
%% \label{thm_main_result}
{\bf Theorem.\,}
\emph{Let $\hat{\mathbf{x}}(\cdot)$ be an estimator for the SLM \eqref{equ_observation_model}, \eqref{equ_SLM_parameter} 
whose mean equals $\bm{\gamma}(\mathbf{x})$ for all $\mathbf{x} \!\in\! \mathcal{X}_{S}$. Then the variance of $\hat{\mathbf{x}}(\cdot)$
at a given parameter vector $\mathbf{x} \!=\! \mathbf{x}_{0} \!\in\! \mathcal{X}_{S}$ 
\vspace{-1mm}
satisfies 
\begin{equation}
\label{equ_lower_vec_est}
 v(\hat{\mathbf{x}}(\cdot);\mathbf{x}_{0}) \ist\geq\ist \sum_{k=1}^{N} L_{\gamma_k}^*\rmv (\mathbf{x}_{0}) \,,
\vspace{-1mm}
\end{equation} 
where $L_{\gamma_k}^*\rmv (\mathbf{x}_{0}) \triangleq \max_{\mathcal{K}_k : |\mathcal{K}_k| = S} L_{\gamma_k}^{\mathcal{K}_k}\rmv (\mathbf{x}_{0})$,
with $L_{\gamma_k}^{\mathcal{K}_k}\rmv (\mathbf{x}_{0})$ given by \eqref{equ_V_CRB_Lfinal} together with \eqref{equ_CRB_LGM_H} and \eqref{equ_s0_expr}.
%% \vspace{2mm}
}
%% \end{theorem}

%% Due to the assumption \eqref{equ_spark_cond}, the component bounds $L_{\gamma_k}^{\mathcal{K}}\rmv (\mathbf{x}_{0})$ and, in turn, 
%% the overall bound $\sum_{k=1}^{N} L_{\gamma_k}^*\rmv (\mathbf{x}_{0})$ are always finite. In certain situations, this fact implies that our bound 
%% will be looser (i.e., smaller) than the CRB derived in \cite{ZvikaCRB}.??? 
%% This occurs, e.g., when there exists no estimator for the SLM having the prescribed bias function $\mathbf{c}(\mathbf{x})$.???

%\mbox{var} \{ g^{k} (\mathbf{y}) \} \geq \max_{|\mathcal{K}| = S, i_{1} =k} \mbox{CRB}^{\mathbf{H}_{\mathcal{K}}}_{\sigma^{2}} \left( \mathbf{x}_{0}^{\mathcal{K}} +  \widetilde{\mathbf{H}}_{\mathcal{K}}\mathbf{H}_{\tiny{supp}(\mathbf{x}_{0}) \setminus \mathcal{K}} \mathbf{x}_{0}^{\tiny{supp}(\mathbf{x}_{0}) \setminus \mathcal{K}} \right)  e^{- \frac{1}{\sigma^{2}} \| (\mathbf{I} - \mathbf{P}_{\mathcal{K}})\mathbf{H}_{\tiny{supp}(\mathbf{x}_{0}) \setminus \mathcal{K}} \mathbf{x}_{0}^{\tiny{supp}(\mathbf{x}_{0}) \setminus \mathcal{K}} \|^{2} }.
%\end{equation}

\vspace{-.8mm}

%%%%%%%%%%%%%%%%%%%%%%%%%%%%%%%%%%%
\section{Special Case: Unbiased Estimation for the SSNM} 
\label{sec_lower_bounds_SSNM}
%%%%%%%%%%%%%%%%%%%%%%%%%%%%%%%%%%%

\vspace{1.3mm}

%% \subsection{Special case: Unbiased estimation for the SSNM}

%\subsubsection{Arbitrary $S$, $N$}
The SSNM in \eqref{equ_observation_model_SSNM} is a special case of the SLM with $\mathbf{H} = \mathbf{I}$.  
We now consider unbiased estimation (i.e., $\bm{\gamma}(\mathbf{x}) \equiv \mathbf{x}$) for the SSNM. 
Since an unbiased estimator with uniformly minimum variance does not exist \cite{AlexZvikaICASSP}, 
we are interested in 
%% the LMVU 
a lower variance bound at a fixed 
%% parameter vector 
$\mathbf{x}_{0} \!\in\! \mathcal{X}_{S}$. 
%% We assume that with ${\|\mathbf{x}_0 \|}_{0} \equiv | \rmv\supp(\mathbf{x}_{0}) | \rmv=\rmv S$. 
We denote by $\xi(\mathbf{x}_{0})$ and $j(\mathbf{x}_{0})$ the value and index, respectively, 
of the $S$-largest (in magnitude) entry of $\mathbf{x}_{0}$; note that this is the  
smallest (in magnitude) nonzero entry of $\mathbf{x}_{0}$ if ${\| \mathbf{x}_{0} \|}_{0} \!=\! S$, and zero if ${\| \mathbf{x}_{0} \|}_{0} \!<\! S$. 

Consider an unbiased estimator $\hat{x}_{k}(\cdot)$.
%%  of the $k\ist$th component $x_k$. 
For $k \in \supp(\mathbf{x}_{0})$, 
using the lower bound $L_{\gamma_k}^{\mathcal{K}_k}\rmv (\mathbf{x}_{0})$ in $\eqref{equ_V_CRB_Lfinal}$
with any index set $\mathcal{K}_k$ of size $|\mathcal{K}_k| \!=\! S$ such that $\supp(\mathbf{x}_{0}) \rmv\rmv\subseteq\rmv \mathcal{K}_k$, 
one can show that 
%% the variance of $\hat{x}_{k}(\cdot)$ satisfies
\begin{equation}
\label{equ_SSNM_bound_k_in}
v(\hat{x}_{k}(\cdot);\mathbf{x}_{0}) \ist\geq\ist \sigma^{2} \ist, \quad\, k \in \supp(\mathbf{x}_{0}) \,.
\end{equation}
This bound is actually the minimum variance (i.e., the variance of the LMVU estimator)
%% \emph{maximally tight} 
since it is achieved by the specific unbiased estimator $\hat{x}_{k} (\mathbf{y}) = y_{k}$ (which is the LMVU estimator for $k \in \supp(\mathbf{x}_{0})$). 
On the other hand, for $k \notin \supp(\mathbf{x}_{0})$, the lower bound $L_{\gamma_k}^{\mathcal{K}_k}\rmv (\mathbf{x}_{0})$ 
with $\mathcal{K}_k = \big( \rmv\rmv\supp(\mathbf{x}_{0}) \!\setminus\! \{j(\mathbf{x}_{0})\} \big) \cup \{ k \}$ can be shown to lead to the 
\vspace{-.3mm}
inequality 
\begin{equation} 
\label{equ_SSNM_bound_k_notin}
v(\hat{x}_{k}(\cdot);\mathbf{x}_{0}) \ist\geq\ist \sigma^{2} \ist e^{-\xi^{2}(\mathbf{x}_{0})/\sigma^{2} }, \quad\, k \notin \supp(\mathbf{x}_{0}) \,.
\vspace{.5mm}
\end{equation}
Combining \eqref{equ_SSNM_bound_k_in} and \eqref{equ_SSNM_bound_k_notin}, a lower bound on the overall variance 
$v(\hat{\x}(\cdot);\mathbf{x}_{0}) = \sum_{k=1}^N v(\hat{x}_{k}(\cdot);\mathbf{x}_0)$ is obtained as
\begin{align}
v(\hat{\x}(\cdot);\mathbf{x}_{0}) &\,\geq \sum_{k \in \supp(\mathbf{x}_{0})} \!\!\sigma^{2} \, +\! \sum_{k \notin \supp(\mathbf{x}_{0})} \!\! \sigma^{2} \ist e^{-\xi^{2}(\mathbf{x}_{0})/\sigma^{2} }.  
%  \rule{9mm}{0mm} % \nonumber\\[.8mm]
		% &\hspace*{15mm}  \eq S \ist \sigma^{2} \rmv+\ist (N\!\rmv-\! S) \ist\ist \sigma^{2} \ist e^{-\xi^{2}(\mathbf{x}_{0})/\sigma^{2} } \rmv\rmv.  
\label{equ_lower_bound_unbiased_ssnm}
\end{align}	
%Because for unbiased estimators the variance equals the MSE, 
%% this is also a lower bound on the MSE at $\mathbf{x}_{0}$. We thus have the following corollary to our theorem of Section \ref{sec:another_lowerbound}.
Thus, recalling that $v(\hat{\x}(\cdot);\mathbf{x}_{0}) = \varepsilon( \hat{\x}(\cdot); \mathbf{x}_{0} )$ for unbiased estimators, we arrive at the following %%% BLIBLABLU
\vspace{2mm}
result.

\noindent 
%% \begin{theorem} 
%% \label{thm_main_result}
{\bf Corollary.\,}
\emph{Let $\hat{\mathbf{x}}(\cdot)$ be an unbiased estimator for the SSNM in \eqref{equ_observation_model_SSNM}. Then the MSE of $\hat{\mathbf{x}}(\cdot)$
at a given 
%% parameter vector 
$\mathbf{x} \!=\! \mathbf{x}_{0} \!\in\! \mathcal{X}_{S}$ satisfies} 
\be
\varepsilon( \hat{\x}(\cdot); \mathbf{x}_{0} ) \,\geq\, \big[ S +\ist (N\!\rmv-\! S) \ist\ist e^{-\xi^{2}(\mathbf{x}_{0})/\sigma^{2} } \big] \ist \sigma^{2} .
\label{equ_lower_bound_unbiased_ssnm_mse}
%% \vspace{-3mm}
\pagebreak %%%%%%%%
\ee	

%where the first component is due to the estimation of the coefficients of $\mathbf{x}$ within the support of $\mathbf{x}_{0}$ and the second term is due to the estimation 
%of the coefficients not in the support of $\mathbf{x}_{0}$. 
This lower bound 
%% \eqref{equ_lower_bound_unbiased_ssnm_mse} 
is tighter (i.e., higher) than the lower bound derived in \cite{AlexZvikaICASSP}. 
Furthermore, in contrast to the bound in \cite{AlexZvikaICASSP}, it is a function of $\mathbf{x}_{0}$ that is everywhere continuous. This fact is theoretically pleasing since the 
MSE of any estimator is a continuous function of $\mathbf{x}_{0}$ \cite{LC}.  %%% BLIBLABLU

%\subsubsection{$S=1$}

Let us consider the special case of 
%% sparsity degree 
$S\!=\! 1$. Here, $\xi(\mathbf{x}_{0})$ and $j(\mathbf{x}_{0})$ are simply the value and index, respectively, of the single nonzero entry of $\mathbf{x}_{0}$.
Using RKHS theory, one can show that the estimator $\hat{\x}(\cdot)$ given componentwise by 
\[ 
\hat{x}_k(\mathbf{y}) \rmv\eq\rmv \begin{cases} y_{j(\mathbf{x}_{0})} \,,  & k =j(\mathbf{x}_{0})\\[.5mm]
  \alpha(\mathbf{y};\mathbf{x}_{0}) \ist\ist y_{k} \,, &  \text{else} \,,
  %% k \in \{1,\ldots,N\} \!\setminus\! \{j(\mathbf{x}_{0})\} \,,
\end{cases}
\]
with $\alpha(\mathbf{y};\mathbf{x}_{0}) \triangleq \exp \rmv\rmv\big( \!\rmv-\rmv\rmv \frac{1}{2 \sigma^{2}} \ist [2 \ist y_{j(\mathbf{x}_{0})} \ist \xi(\mathbf{x}_{0}) + \xi^{2}(\mathbf{x}_{0}) ] \big)$,
\vspace{.4mm}
is the LMVU estimator at $\mathbf{x}_{0}$. That is, the estimator $\hat{\x}(\cdot)$ is unbiased and its MSE achieves the lower bound \eqref{equ_lower_bound_unbiased_ssnm_mse}.
This also means that \eqref{equ_lower_bound_unbiased_ssnm_mse} is actually the minimum MSE (achieved by the LMVU estimator). 
While $\hat{\x}(\cdot)$ is not very practical since it explicitly involves the unknown true parameter $\mathbf{x}_{0}$, its existence demonstrates the tightness of the bound 
\eqref{equ_lower_bound_unbiased_ssnm_mse}.
%% in the present setting. %BLIBLABLU
%Of course, this estimator is not very practical because it explicitly involves the unknown true parameter $\mathbf{x}_{0}$.
%% Further note that for the case $S\!=\! 1$ considered,  
%% $\xi(\mathbf{x}_{0})$ and $j(\mathbf{x}_{0})$ are the value and index, respectively, of the single nonzero component of $\mathbf{x}_{0}$.    

\vspace{-.5mm}

%% \newpage %%%%%%%%%%%

%%%%%%%%%%%%%%%%%%%%%%%%%%%%%%%%%%%%%%%%%%%%%%%%%%%%%%%%%%
\section{Numerical Results} 
\label{sec_numerical}
%%%%%%%%%%%%%%%%%%%%%%%%%%%%%%%%%%%%%%%%%%%%%%%%%%%%%%%%%%
%\subsection{Comparsion with the Thresholding- and ML-estimators} 

\vspace{1mm}

For the 
%% special case of the 
SSNM in \eqref{equ_observation_model_SSNM},
%%  with $S \!=\! 1$, 
we will compute the lower variance 
bound $\sum_{k=1}^{N} \rmv\rmv L_{\gamma_k}^*\rmv (\mathbf{x}_{0})$
(see \eqref{equ_lower_vec_est}) and compare it with the 
%% actual 
variance of two established estimators, namely, the maximum likelihood (ML) estimator and the hard-thresholding (HT) estimator.
The ML estimator is given by 
\[
%% \label{equ_ML_est}
\ML_est(\mathbf{y}) \,\triangleq\, \argmax_{\mathbf{x}' \in \mathcal{X}_{S}} f ( \mathbf{y};  \mathbf{x}' ) \ist=\, {\mathsf P}_{\!S} ( \mathbf{y}  ) \,,
\]
where the operator $\mathsf{P}_{\! S}$ retains the $S$ largest (in magnitude) entries and zeros out
all others. The HT estimator $\hat{\mathbf{x}}_{\text{HT}}  (\mathbf{y})$ is given 
\vspace{-.5mm}
by
\begin{equation}
\label{equ_def_thr_func}
\hat{x}_{\text{HT},k}  (\mathbf{y}) \eq \begin{cases} y_k \,,  & |y_k| \geq T\\[.5mm]
  0 \,, & \text{else} \ist,
\end{cases} 
\end{equation} 
where $T$ is a fixed threshold.
%% A popular choice for the threshold is $T= \sigma \sqrt{2 \log_{e} N}$  \cite{Mallat2008}. 

For simplicity, we consider the SSNM for $S \!=\! 1$.
In this case, the bound \eqref{equ_lower_vec_est} can be shown to 
\vspace{.5mm}
be
\begin{equation}
\label{equ_lower_bound_ml_ht}
v(\hat{\mathbf{x}}(\cdot); \mathbf{x}_0) \ist\geq\ist\ist L_{\gamma_{j}}^{\mathcal{K}_j}\rmv (\mathbf{x}_{0})
\ist+\ist (N \!-\! 1) \, e^{ - \xi^2(\mathbf{x}_{0}) / \sigma^{2}} L_{\gamma_{i}}^{\mathcal{K}_i}\rmv (\mathbf{x}_{0}) \,, 
\vspace{.5mm}
\end{equation} 
where $j \!\triangleq\! j(\mathbf{x}_{0})$, $i$ is any index different from $j(\mathbf{x}_{0})$ (it can be shown that all such indices equally maximize
the lower bound), $\mathcal{K}_{j} \!\triangleq\! \{ j(\mathbf{x}_{0}) \}$, and $\mathcal{K}_{i} \!\triangleq\! \{i \}$.
(We note that \eqref{equ_lower_bound_ml_ht} simplifies to \eqref{equ_lower_bound_unbiased_ssnm} for the special case of an unbiased estimator.)
%% As always, the function $\bm{\gamma}(\mathbf{x})$ involved in the bound 
%% %% on the right-hand side 
%% is the mean of the estimator $\hat{\mathbf{x}}(\cdot)$. 
%%   lower variance bound holds for any estimator $\hat{\mathbf{x}}(\cdot)$ with mean $\bm{\gamma}(\mathbf{x})$.
Since we 
%% consider specifically 
compare the bound \eqref{equ_lower_bound_ml_ht} to the ML and HT estimators, $\bm{\gamma}(\mathbf{x})$ is set equal to the 
mean of the respective estimator (ML or HT).

\begin{figure}
\vspace{-1mm}
\centering
\psfrag{SNR}[c][c][.9]{\uput{3.4mm}[270]{0}{\hspace{3mm}SNR [dB]}}
\psfrag{title}[c][c][.9]{\uput{2.5mm}[270]{0}{}}
\psfrag{x_0}[c][c][.9]{\uput{0.3mm}[270]{0}{$0$}}
\psfrag{x_0_001}[c][c][.9]{\uput{0.3mm}[270]{0}{$-30$}}
\psfrag{x_0_01}[c][c][.9]{\uput{0.3mm}[270]{0}{$-20$}}
\psfrag{x_0_1}[c][c][.9]{\uput{0.3mm}[270]{0}{$-10$}}
\psfrag{x_1}[c][c][.9]{\uput{0.3mm}[270]{0}{$0$}}
\psfrag{x_10}[c][c][.9]{\uput{0.3mm}[270]{0}{$10$}}
\psfrag{x_100}[c][c][.9]{\uput{0.3mm}[270]{0}{$20$}}
\psfrag{y_0}[c][c][.9]{\uput{0.1mm}[180]{0}{$0$}}
\psfrag{y_1}[c][c][.9]{\uput{0.1mm}[180]{0}{}}
\psfrag{y_10}[c][c][.9]{\uput{0.1mm}[180]{0}{$10$}}
\psfrag{y_2}[c][c][.9]{\uput{0.1mm}[180]{0}{$2$}}
\psfrag{y_3}[c][c][.9]{\uput{0.1mm}[180]{0}{}}
\psfrag{y_4}[c][c][.9]{\uput{0.1mm}[180]{0}{$4$}}
\psfrag{y_5}[c][c][.9]{\uput{0.1mm}[180]{0}{}}
\psfrag{y_6}[c][c][.9]{\uput{0.1mm}[180]{0}{$6$}}
\psfrag{y_7}[c][c][.9]{\uput{0.1mm}[180]{0}{}}
\psfrag{y_8}[c][c][.9]{\uput{0.1mm}[180]{0}{$8$}}
\psfrag{y_9}[c][c][.9]{\uput{0.1mm}[180]{0}{}}
\psfrag{variancebound}[c][c][.9]{\uput{1mm}[90]{0}{variance/bound}}
\psfrag{bML}[l][l][0.8]{bound on $v(\ML_est(\cdot);\mathbf{x}_0)$}
\psfrag{ML}[l][l][0.8]{$v(\ML_est(\cdot);\mathbf{x}_0)$}
\psfrag{HT3}[l][l][0.8]{$v(\hat{\mathbf{x}}_{\text{HT}}(\cdot); \mathbf{x}_0)$, $T \!\rmv=\! 3$}
\psfrag{HT4}[l][l][0.8]{$v(\hat{\mathbf{x}}_{\text{HT}}(\cdot); \mathbf{x}_0)$, $T \!\rmv=\! 4$}
\psfrag{HT5}[l][l][0.8]{$v(\hat{\mathbf{x}}_{\text{HT}}(\cdot); \mathbf{x}_0)$, $T \!\rmv=\! 5$}
\psfrag{bHT3}[l][l][0.8]{bound on $v(\hat{\mathbf{x}}_{\text{HT}}(\cdot); \mathbf{x}_0)$, $T \!\rmv=\! 3$}
\psfrag{bHT4}[l][l][0.8]{bound on $v(\hat{\mathbf{x}}_{\text{HT}}(\cdot); \mathbf{x}_0)$, $T \!\rmv=\! 4$}
\psfrag{bHT5}[l][l][0.8]{bound on $v(\hat{\mathbf{x}}_{\text{HT}}(\cdot); \mathbf{x}_0)$, $T \!\rmv=\! 5$}
\psfrag{MLarrow}[l][l][0.8]{ML}
\psfrag{T3}[l][l][0.8]{HT \!($T \!\rmv=\! 3$)}
\psfrag{T4}[l][l][0.8]{HT \!($T \!\rmv=\! 4$)}
\psfrag{T5}[l][l][0.8]{HT \!($T \!\rmv=\! 5$)}
\centering
\hspace*{-0mm}\includegraphics[height=5.2cm,width=9cm]{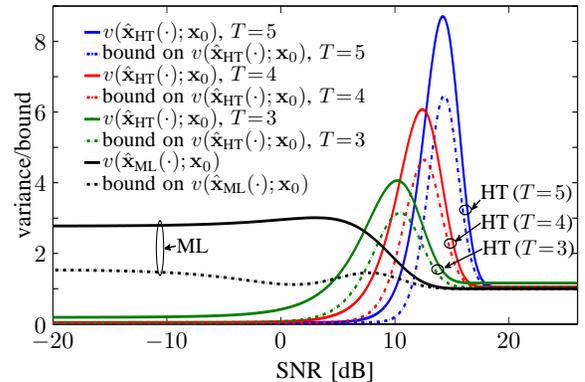}
\vspace{-2.5mm}
  \caption{Variance of the ML and HT estimators and corresponding lower bounds versus the SNR $\xi^{2}(\mathbf{x}_{0})/\sigma^{2}\rmv\rmv$,
  for the SSNM
  %% SLM with $\mathbf{H}=\mathbf{I}$, 
  with $N\!=\!5$ and $S\!=\!1$.} %%  and $\sigma^{2}=1$
\label{fig_bounds_1}
\vspace*{-2mm}
\end{figure}

For a numerical evaluation, we generated parameter vectors $\mathbf{x}_0$ with $N\!\rmv=\!5$, $S\!=\!1$, $j(\mathbf{x}_{0}) \!=\! 1$, and different $\xi(\mathbf{x}_{0})$.
(The fixed choice $j(\mathbf{x}_{0}) \!=\! 1$ is justified by the fact that neither the variances of the ML and HT estimators nor the corresponding variance bounds depend on $j(\mathbf{x}_{0})$.)
In Fig.~\ref{fig_bounds_1}, we plot the variances $v(\ML_est(\cdot); \mathbf{x}_0)$ and $v( \hat{\mathbf{x}}_{\text{HT}} (\cdot); \mathbf{x}_0)$ 
(the latter for three different choices of 
%% the threshold parameter 
$T$ in \eqref{equ_def_thr_func}) along with the corresponding bounds \eqref{equ_lower_bound_ml_ht}, as a function of the signal-to-noise ratio (SNR) 
$\xi^{2}(\mathbf{x}_{0})/\sigma^{2}$. It is seen that for SNR larger than about 18 dB, %%% BLIBLABLU
%% 15\,dB, 
all variances and bounds are effectively equal (for the HT estimator, 
this is true if $T$ is not too small). However, in the medium-SNR
%% intermediate SNR 
range, the variances of the ML and HT estimators are significantly higher than
%% there is a significant gap between the actual variance performance of the ML and HT estimators and 
the corresponding lower bounds. 
We can conclude that there \emph{might} exist estimators with the same mean as that of the ML or HT estimator but smaller variance.
%% (or MSE).  
Note, however, that a positive statement regarding the existence of such estimators cannot be based on our analysis.

\vspace{-.5mm}

%% \newpage %%%%%%%%%

%%%%%%%%%%%%%%%%%%%%%%%%%%%%%%%%%%%%%%%%%%%%%%%%%%%%%%%%%%
\section{Conclusion} 
%%%%%%%%%%%%%%%%%%%%%%%%%%%%%%%%%%%%%%%%%%%%%%%%%%%%%%%%%%

\vspace{.8mm}

Using the mathematical framework of reproducing kernel Hilbert spaces,
%%  (RKHS), 
we derived a novel lower bound on the variance of estimators 
of a sparse vector under a bias constraint. The observed vector was assumed to be a linearly transformed and noisy version of the sparse vector to be estimated.
This setup includes the underdetermined case relevant to compressed sensing. In the special case of unbiased estimation of a noise-corrupted sparse vector,
our bound improves on the best known lower bound. A comparison 
%% of our bound 
with the variance of two established estimators 
%% (maximum likelihood and hard thresholding)
showed that there might exist estimators with the same bias 
%% as these estimators 
but a smaller variance.

\vspace{-1mm}

%% \pagebreak %%%%%%%%%

\bibliographystyle{ieeetr}
%% \bibliographystyle{plain}
%\bibliography{LitAJPHD}
\bibliography{/Users/ajung/Studium/Diss/PHDThesis/LitAJPHD.bib}

\end{document}